\documentclass[reqno,12pt]{amsart}
\usepackage{ifthen}
\usepackage[pdftex]{graphicx}
\DeclareGraphicsExtensions{.pdf}

\newcommand{\version}{journal}

\ifthenelse{\equal{\version}{journal}}%
{\message{
--- journal version ---
}}%
{\message{
--- personal version ---
}
}

\usepackage{amssymb}
\usepackage{amsmath}
\usepackage{amsfonts}
\usepackage{latexsym}

\numberwithin{equation}{section}
\numberwithin{table}{section}
\numberwithin{figure}{section}
\newtheorem{theorem}{Theorem}[section]

\newtheorem{remark}[theorem]{Remark}
\newtheorem{example}[theorem]{Example}

\newtheorem{algorithm}[theorem]{Algorithm}

\newcommand{\boldc}{{\mathbf c}}

\newcommand{\bbR}{{\mathbb{R}}}
\newcommand{\R}{{\mathbb{R}}}

\newcommand{\Z}{{\mathbb{Z}}}
\newcommand{\boldm}{{\boldsymbol{m}}}
\newcommand{\boldn}{{\boldsymbol{n}}}

\newcommand{\boldPhi}{{\boldsymbol{\Phi}}}
\newcommand{\boldphi}{{\boldsymbol{\phi}}}
\newcommand{\boldPsi}{{\boldsymbol{\Psi}}}
\newcommand{\boldpsi}{{\boldsymbol{\psi}}}

\newcommand{\<}{\langle}
\renewcommand{\>}{\rangle}
\newcommand{\binomial}[2]{\begin{pmatrix}#1\\#2\end{pmatrix}} %use \dbinom{}{}
\def\dfrac{\displaystyle\frac}

\begin{document}

\title[Two-direction multiwavelet moments]
%[Moments of orthogonal two-direction multiscaling functions]
%{Moments of two-direction multiscaling functions}
%{Moments of two-direction multiscaling functions for dilation factor $d$}
%{Moments of orthogonal two-direction multiscaling functions for dilation factor $d$}
%{Moments of two-direction (bi-)orthogonal multiscaling functions and two-direction multiwavelets for dilation factor $d$}
%{Two-direction multiwavelet moments for dilation factor $d$}
{Two-direction multiwavelet moments}

\author{Soon-Geol Kwon}

\thanks{This paper was supported by (in part) Sunchon National University Research Fund in 2010.}
%\thanks{Research was supported by a nondirected research fund from Sunchon National %University, Korea.}
%\thanks{{\em E-mail address:} sgkwon@sunchon.ac.kr}

\address{Soon-Geol Kwon \\
Department of Mathematics Education, Sunchon National University,
Sunchon, 540-742, Korea}
\email{sgkwon@sunchon.ac.kr}
%\acknowledge{Thanks. }

\begin{abstract}
Two-direction multiscaling functions $\boldphi$ and two-direction multiwavelets
$\boldpsi$ associated with $\boldphi$ are a more general and more flexible setting
than one-direction multiscaling functions and multiwavelets.
% associated with multiscaling functions.
%
In this paper, we derive two methods for computing continuous moments
of orthogonal two-direction multiscaling functions $\boldphi$ and
orthogonal two-direction multiwavelets $\boldpsi$ associated with $\boldphi$.
The first method is by doubling and the second method is by separation.
%
%Finally, an
Two examples for both methods are given.
%We derive an algorithm for computing the continuous moments
%of the two-direction multiscaling function $\boldphi$ and
%the corresponding two-direction multiwavelets $\boldpsi$.
\end{abstract}

%\subjclass[2000]{42C40, 42C15, 94A12}
\subjclass[2010]{42C15}

\keywords{moments, continuous moments, two-direction multiscaling functions, two-direction multiwavelets.}

%\date{\today.}
%\date{\today and, in revised form, .}

\maketitle

%%%%%%%%%%%%%%%%%%%%%%%%%%%%%%%
%           main text         %
%%%%%%%%%%%%%%%%%%%%%%%%%%%%%%%

%\bigskip
\section{Introduction}\label{introduction}

Two-direction multiscaling functions $\boldphi$ and two-direction multiwavelets
$\boldpsi$ associated with $\boldphi$, which are a more general setting than
the one-direction multiscaling functions and multiwavelets, are investigated in
~\cite{Kwon2011c,Kwon2011b,WangZhouWang2011,XieYang2006,Yang2006,YangLi2007AMC,YangLi2007SciChina,YangXie2008}.
%~\cite{Kwon2011c,Kwon2011b,WangZhouWang2011,Xie2008,Xie2010,XieYang2006,Yang2006,YangLi2007AMC,YangLi2007SciChina,YangXie2008}.
%~\cite{Xie2010,XieYang2006,YangLi2007AMC,YangXie2008}.
The two-direction setting is more flexible than the one-direction setting.

In multiwavelet theory, computation of continuous moments for two-direction
multiwavelets associated with two-direction multiscaling functions is important,
since vanishing continuous moments for two-direction multiwavelets provide
the approximation order for two-direction multiscaling functions.
The discrete or continuous moments for the orthogonal two-direction multiscaling
functions and multiwavelets are necessary in theory such as for establishing
Condition E and approximation order, and many applications such as solutions of
differential equations, signal processing, and image processing, especially
for prefiltering or balancing two-direction multiscaling functions.
%
%This is a great advantage both in terms of computational cost and
%quality of results in many applications.
%
It is well-known that continuous moments of one-direction orthogonal
multiscaling functions and multiwavelets can be computed
if their recurrence coefficients are given, for example,
see ~\cite{Gopinath1992,Johnson2000,Keinert2004, Kwon2009}.
%(for e.g., See ~\cite{Gopinath1992,Johnson2000,Keinert2004, Kwon2009})
%~\cite{Keinert2004, Kwon2009})

The main objective of this paper is to derive two different methods for
computing continuous moments of the orthogonal two-direction multiscaling
functions $\boldphi$ and multiwavelets $\boldpsi$ associated with $\boldphi$.
To derive the method for computing continuous moments by doubling,
we investigate the following: for $s=1,2,\ldots,d-1$, where $d$ is a dilation factor
\begin{itemize}
\item investigate recursion formulas for continuous moments of doubled $\boldPhi(x)=[\boldphi(x),\boldphi(-x)]^T$ and $\boldPsi^{(s)}(x)=[\boldpsi^{(s)}(x),\boldpsi^{(s)}(-x)]^T$;
    %for $s=1,2,\ldots,d-1$;
\item investigate normalization for the zeroth continuous moment of $\boldPhi$;
\item investigate the continuous moments of $\boldPhi$ and $\boldPsi^{(s)}$;
\item choose the continuous moments of $\boldphi$ and $\boldpsi^{(s)}$ from the continuous moments of $\boldPhi$ and $\boldPsi^{(s)}$, respectively.
\end{itemize}
To derive the method for computing continuous moments by separation,
we investigate the following: for $s=1,2,\ldots,d-1$, where $d$ is a dilation factor
\begin{itemize}
\item investigate recursion formulas for the continuous moments of $\boldphi$ and $\boldpsi^{(s)}$ by separation;
\item investigate normalization for the zeroth continuous moment of $\boldphi$;
\item investigate the continuous moments of $\boldphi$ and $\boldpsi^{(s)}$.
\end{itemize}
We provide two algorithms, one by doubling and the other by separation, for
computing continuous moments of the orthogonal two-direction multiscaling functions
$\boldphi$ and multiwavelets $\boldpsi$ associated with $\boldphi$.

This paper is organized as follows.
Orthogonal two-direction multiscaling functions $\boldphi$ and multiwavelets
$\boldpsi$ are introduced in section 2. The main results for computing continuous
moments of the orthogonal two-direction multiscaling functions $\boldphi$
and multiwavelets $\boldpsi$ are stated in sections 3 and 4.
In section 3, one method for computing moments by doubling is introduced.
In section 4, the other method for computing moments by separation is introduced.
Finally, two examples for illustrating the general theory in sections 2, 3, and 4
are given in section 5.

\bigskip
\section{Two-direction multiscaling functions and multiwavelets}\label{2direc}

In this section we review orthogonal two-direction multiscaling functions and
multiwavelets associated with the orthogonal two-direction multiscaling functions
(see ~\cite{XieYang2006,Yang2006,YangLi2007AMC,YangLi2007SciChina,YangXie2008}).
%(See ~\cite{Xie2008,Xie2010,XieYang2006,YangLi2007AMC,YangLi2007SciChina,YangXie2008})
%~\cite{Xie2010,YangXie2008,YangLi2007AMC,YangLi2007SciChina}.

A {\em two-direction  multiscaling function} of multiplicity $r$ and dilation
factor $d$ is a vector of $r$ real or complex-valued functions
\begin{equation*}
\boldphi(x) = \left[ \phi_{1}(x), \phi_{2}(x), \ldots, \phi_{r}(x)
\right]^T, \quad x \in \bbR,
\end{equation*}
which satisfies a recursion relation
\begin{equation}\label{recrel}
  \boldphi(x) = \sqrt{d} \,\sum_{k\in\Z} \left[ P^{+}_{k} \,\boldphi(dx-k) + P^{-}_{k} \,\boldphi(k-dx) \right]
\end{equation}
and generates a multiresolution approximation of $L^{2}(\bbR)$.
The $P^{+}_{k}$, $P^{-}_{k}$, called positive- and negative-direction
recursion coefficients for $\boldphi$, respectively, are $r \times r$ matrices.

{\em Two-direction multiwavelets} $\boldpsi^{(s)}$, $s=1,2,\ldots,d-1$, associated with $\boldphi$
satisfy
\begin{equation}\label{recrelpsi}
  \boldpsi^{(s)}(x) = \sqrt{d} \,\sum_{k\in\Z} \left[ Q^{(s)+}_{k} \,\boldphi(dx-k) + Q^{(s)-}_{k} \,\boldphi(k-dx) \right].
\end{equation}
The $Q^{(s)+}_{k}$, $Q^{(s)-}_{k}$ are called positive- and negative-direction
recursion coefficients for $\boldpsi^{(s)}$, $s=1,2,\ldots,d-1$, respectively.
%, are $r \times r$ matrices.

Two-direction multiscaling function $\boldphi$ and multiwavelets $\boldpsi^{(s)}$
associated with $\boldphi$ are {\em orthogonal} if for all $j,k\in\Z$, and
$s,t=1,2,\ldots,d-1$,
\begin{align*}
  \< \boldphi(x-j), \boldphi(x-k) \> &= \delta_{jk} I_{r}, \\
  \< \boldphi(x-j), \boldphi(k-x) \> &= \bold{O}, \\
  \< \boldpsi^{(s)}(x-j), \boldpsi^{(t)}(x-k) \> &= \delta_{st} \delta_{jk} I_{r}, \\
  \< \boldpsi^{(s)}(x-j), \boldpsi^{(t)}(k-x) \> &= \bold{O}, \\
  \< \boldphi(x-j), \boldpsi^{(s)}(x-k) \> &= \bold{O}, \\
  \< \boldphi(x-j), \boldpsi^{(s)}(k-x) \> &= \bold{O},
\end{align*}
where $I_{r}$ is the $r \times r$ identity matrix and $\bold{O}$ is the $r\times r$
zero matrix. Here the inner product is defined by
\begin{equation*}
  \< \boldphi, \boldpsi^{(s)} \> = \int \boldphi(x) {\boldpsi^{(s)}}^*(x) \;\textrm{d}x,
\end{equation*}
where * denotes the complex conjugate transpose.
This inner product is an $r \times r$ matrix.

By taking the Fourier transform on both sides of ~\eqref{recrel}, we have
\begin{equation}\label{ft_recrel}
    \widehat{\boldphi}(\xi) = P^{+}(e^{-i\xi/d}) \,\widehat{\boldphi}(\xi/d) + P^{-}(e^{-i\xi/d}) \,\overline{\widehat{\boldphi}(\xi/d)},
\end{equation}
where
\begin{equation}
    P^{+}(z) = \frac{1}{\sqrt{d}} \,\sum_{k\in\Z} P^{+}_k z^k
    \qquad \text{and} \qquad
    P^{-}(z) = \frac{1}{\sqrt{d}} \,\sum_{k\in\Z} P^{-}_k z^k
\end{equation}
are called positive- and negative-direction mask symbols, respectively.
%Then the refinement mask for $\boldphi$ is
%\begin{equation}\label{Hphi}
%    \bold{P}_{\boldphi}(z) = \begin{bmatrix}
%                               P^{+}(z) & P^{-}(z) \\
%                             \end{bmatrix}.
%\end{equation}

%In order to investigate the existence of the solutions of the two-direction
%matrix refinable equation ~\eqref{recrel},
We rewrite the two-direction recursion relation ~\eqref{recrel} as
\begin{equation}\label{recrelneg}
  \boldphi(-x) = \sqrt{d} \,\sum_{k\in\Z} \left[ P^{+}_{k} \,\boldphi(-dx-k) + P^{-}_{k} \,\boldphi(k+dx) \right]
\end{equation}

By taking the Fourier transform on both sides of ~\eqref{recrelneg}, we have
\begin{equation}\label{ft_recrelneg}
    \overline{\widehat{\boldphi}(\xi)} = \overline{P^{+}(e^{-i\xi/d})}\, \overline{\widehat{\boldphi}(\xi/d)} + \overline{P^{-}(e^{-i\xi/d})} \, \widehat{\boldphi}(\xi/d).
\end{equation}
From ~\eqref{ft_recrel} and ~\eqref{ft_recrelneg}, we have
\begin{equation}\label{Phifreq}
    \begin{bmatrix}
      \widehat{\boldphi}(\xi) \\ \noalign{\medskip}
      \overline{\widehat{\boldphi}(\xi)} \\
    \end{bmatrix}
    =
    \begin{bmatrix}
      P^{+}(e^{-i\xi/d}) & P^{-}(e^{-i\xi/d})  \\ \noalign{\medskip}
      \overline{P^{-}(e^{-i\xi/d})} & \overline{P^{+}(e^{-i\xi/d})} \\
    \end{bmatrix}
    \begin{bmatrix}
      \widehat{\boldphi}(\xi/d) \\ \noalign{\medskip}
      \overline{\widehat{\boldphi}(\xi/d)} \\
    \end{bmatrix}.
\end{equation}
We see that ~\eqref{recrel} has a solution if and only if ~\eqref{Phifreq}
has a solution. That is, ~\eqref{recrel} is equivalent to ~\eqref{Phifreq}.

Let \begin{equation}\label{recrelPhi}
\boldPhi(x) =
\begin{bmatrix}
  \boldphi(x) \\ \noalign{\smallskip}
  \boldphi(-x) \\
\end{bmatrix}
= \sqrt{d} \sum_{k\in\Z}
\begin{bmatrix}
    P^{+}_k & P^{-}_k \\ \noalign{\medskip}
    P^{-}_{-k} & P^{+}_{-k} \\
\end{bmatrix}
\boldPhi(dx-k)
\end{equation}
be the matrix refinable equation of $\boldPhi$.
Equation ~\eqref{recrelPhi} is called the {\em deduced} $d$-scale matrix refinement
equation of the two-direction refinement equation ~\eqref{recrel}.
%Its solution $\boldPhi$ of ~\eqref{recrelPhi} is called the {\em deduced} $d$-scale
%matrix refinement equation of $\boldPhi$.
Then ~\eqref{Phifreq} is the $d$-scale matrix refinement equation in the frequency
domain of $\boldPhi$. Its refinement mask is
\begin{equation}\label{refinemask}
\bold{P}_{\boldPhi}(z) = \begin{bmatrix}
      P^{+}(z) & P^{-}(z) \\ \noalign{\medskip}
      \overline{P^{-}(z)} & \overline{P^{+}(z)} \\
    \end{bmatrix}.
\end{equation}
%for $|z|=1$ and $z\in\C$.
If $P^{+}_{k}$ and $P^{-}_{k}$, $k\in\Z$, are real, then
\begin{equation}\label{refinemaskreal}
\bold{P}_{\boldPhi}(z) = \begin{bmatrix}
      P^{+}(z) & P^{-}(z) \\ \noalign{\medskip}
      \overline{P^{-}(z)} & \overline{P^{+}(z)} \\
    \end{bmatrix}
= \frac{1}{\sqrt{d}} \sum_{k\in\Z}
\begin{bmatrix}
    P^{+}_k & P^{-}_k \\ \noalign{\medskip}
    P^{-}_{-k} & P^{+}_{-k} \\
\end{bmatrix} z^k.
\end{equation}
In this paper we only consider real recursion coefficients $P^{+}_{k}$, $P^{-}_{k}$, $Q^{(s)+}_{k}$, and $Q^{(s)-}_{k}$ in $\R^{r\times r}$ for $s=1,2,\ldots,d-1$ and $k\in\Z$.

{\bf Condition E.}
A matrix $A$ satisfies {\em Condition E} if it has a simple eigenvalue of 1,
and all other eigenvalues are smaller than 1 in absolute value.

Condition E for $\boldphi$ defined in ~\eqref{recrel} means that the matrix
\begin{equation}\label{condsymbol}
\begin{bmatrix}
        P^{+}(1) & P^{-}(1) \\
        P^{-}(1) & P^{+}(1) \\
    \end{bmatrix}
= \frac{1}{\sqrt{d}} \sum_{k\in\Z}
    \begin{bmatrix}
        P^{+}_k & P^{-}_k \\ \noalign{\medskip}
        P^{-}_{k} & P^{+}_{k} \\
    \end{bmatrix}
\end{equation}
%\begin{align}\label{condsymbol}
%    P^{+}(1) + P^{-}(1) = \frac{1}{\sqrt{d}} \sum_{k\in\Z} [ P^{+}_k + P^{-}_k ]
%\end{align}
satisfies Condition E (see ~\cite[Theorem 3]{XieYang2006} for $d=2$). %(see ~\cite{Kwon2011c}).

For the scalar case $r=1$, since the eigenvalues of the matrix in ~\eqref{condsymbol} are
$P^{+}(1)+P^{-}(1)$ and $P^{+}(1)-P^{-}(1)$, Condition E for $\boldphi$
is equivalent to
\begin{equation*}%\label{conditionE1}
  P^{+}(1) + P^{-}(1) = 1, \qquad |P^{+}(1) - P^{-}(1)| < 1,
\end{equation*}
or
\begin{equation*}%\label{conditionE2}
  P^{+}(1) - P^{-}(1) = 1, \qquad |P^{+}(1) + P^{-}(1)| < 1,
\end{equation*}
that is, %for the scalar case $r=1$ Condition E for $\boldphi$ is equivalent to
\begin{equation*}%\label{conditionE3}
  \frac{1}{\sqrt{d}} \sum_{k\in\Z} [ P^{+}_k + P^{-}_k ] = 1, \qquad
  \frac{1}{\sqrt{d}} \left| \sum_{k\in\Z} [ P^{+}_k - P^{-}_k ] \right| < 1,
\end{equation*}
or
\begin{equation*}
  \frac{1}{\sqrt{d}} \sum_{k\in\Z} [ P^{+}_k - P^{-}_k ] = 1, \qquad
  \frac{1}{\sqrt{d}} \left| \sum_{k\in\Z} [ P^{+}_k + P^{-}_k ] \right| < 1
\end{equation*}
(see ~\cite[Theorem 2]{YangLi2007AMC}).

%For the scalar case $r=1$,
%\begin{equation*}
%    P^{+}(1) + P^{-}(1) = \frac{1}{\sqrt{d}} \sum_{k\in\Z} [ P^{+}_k + P^{-}_k ] = 1.
%\end{equation*}

It is well-known that if $\boldphi$ is a compactly supported $L^2$-stable solution
of ~\eqref{recrel}, then Condition E for $\boldphi$ is satisfied
(see ~\cite[Theorem 3]{XieYang2006} for $d=2$).

Condition E for $\boldphi$ is the condition for the stability of
the multiresolution approximation produced by $\boldphi$ or for the existence
of $\boldphi$ in ~\eqref{recrel} (see ~\cite[Theorem 3]{Plonka1998}).

{\bf Approximation Order.}
The two-direction multiscaling function approximation to a function $f$ at
resolution $d^{-j}$ is given by the series
\begin{equation}\label{fseries}
  P_{j}f = \sum_{k\in\Z} [ \<f, \boldphi^{+}_{jk} \> \,\boldphi^{+}_{jk}
  + \<f, \boldphi^{-}_{jk} \> \,\boldphi^{-}_{jk} ],
\end{equation}
where
\begin{gather*}
  \boldphi^{+}_{jk}(x) = d^{-j/2}\,\boldphi(d^{j}x - k)
  \qquad  \text{and} \qquad
  \boldphi^{-}_{jk}(x) = d^{-j/2}\,\boldphi(k - d^{j}x).
\end{gather*}
The two-direction multiscaling function $\boldphi$ provides
{\em approximation order $p$} if
\begin{equation*}
    \| f - P_{n}f \| = O(d^{-np}),
\end{equation*}
whenever $f$ has $p$ continuous derivatives.

The two-direction multiscaling function $\boldphi$ has {\em accuracy $p$}
if all polynomials of degree up to $p-1$ can be expressed locally as linear
combinations of integer shifts of $\boldphi(x)$ and $\boldphi(-x)$. That is,
there exist row vectors $(\boldc^{+}_{j,k})^{*}$ and $(\boldc^{-}_{j,k})^{*}$,
$j=0,\ldots{},p-1$, $k \in \Z$, so that
\begin{equation}\label{approxint}
  x^{j} = \sum_{k\in\Z} [(\boldc^{+}_{j,k})^{*} \boldphi(x-k)
  + (\boldc^{-}_{j,k})^{*} \boldphi(k-x)].
\end{equation}
%The recursion coefficients $\{P^{+}_k, P^{-}_k\}$ of a matrix refinement
%equation satisfy the {\em sum rules of order p} if there exist vectors
%$\boldy_{0},\boldy_{1}\ldots,\boldy_{p-1}$ with $\boldy_{0}\neq 0$, which satisfy

%According to multiwavelet theory,
%As shown in ~\cite{Jia1997},
It is well-known that $\boldphi$ provides approximation order $p$
if and only if $\boldphi$ has accuracy $p$ (see~\cite{Keinert2004,Plonka1998}).

A high approximation order is desirable in applications. A minimum
approximation order of 1 is a required condition in many theorems.

Throughout this paper we assume that the two-direction multiscaling function
$\boldphi$ is orthogonal, has compact support, is continuous (which implies
approximation order at least 1), and satisfies Condition E.

%%%%%%%%%%%%%%%%%%%%%%%%%%%%%%%%%%%%%%%%%%%%%%%%%%%%%%%%%%%%%%%%%%%%%%%%%%%%%%%%%%%
\bigskip
\section{Moments of two-direction multiwavelets by doubling}\label{moments}
%%%%%%%%%%%%%%%%%%%%%%%%%%%%%%%%%%%%%%%%%%%%%%%%%%%%%%%%%%%%%%%%%%%%%%%%%%%%%%%%%%%

It is well-known that the discrete and continuous moments of one-direction
orthogonal (or biorthogonal) multiscaling functions and multiwavelets associated
with multiscaling functions can be computed if the recurrence coefficients of
multiscaling functions and multiwavelets are given, for example,
see ~\cite{Gopinath1992,Johnson2000,Keinert2004, Kwon2009}.

%In this section we investigate how to compute continuous moments of the
%two-direction multiscaling function $\boldphi$ and the corresponding
%two-direction multiwavelet $\boldpsi$ can be computed if the recurrence
%coefficients $P^{+}_k$, $P^{-}_k$ of $\boldphi$ and $Q^{(s)+}_k$, $Q^{(s)-}_k$ of %$\boldpsi$ are given.

The discrete or continuous moments for the orthogonal two-direction
multiscaling functions and multiwavelets are necessary in theory such as
for establishing Condition E, approximation order, {\em etc}.

In this section we investigate how to compute the integral $\boldm_{0}$, the zeroth
continuous moment, and higher moments of the orthogonal two-direction multiscaling
functions $\boldphi$ in terms of the recurrence coefficients $P^{+}_k$, $P^{-}_k$
of $\boldphi$.

Similarly, we investigate how to compute the continuous moments of orthogonal two-direction multiwavelets $\boldpsi^{(s)}$ associated with the orthogonal two-direction
multiscaling function $\boldphi$ in terms of the recurrence coefficients
$Q^{(s)+}_k$, $Q^{(s)-}_k$ of $\boldpsi^{(s)}$ for $s=1,2,\ldots,d-1$.

%We begin by defining some terms.
%We recall from ~\eqref{refinemaskreal} that the symbol of the two-direction
%multiscaling function $\boldPhi(x)$ is
%\begin{equation*}
%\bold{P}_{\boldPhi}(z) = \frac{1}{\sqrt{d}} \sum_{k\in\Z}
%\begin{bmatrix}
%    P^{+}_k & P^{-}_k \\ \noalign{\smallskip}
%    P^{-}_{-k} & P^{+}_{-k} \\
%\end{bmatrix} z^k.
%\end{equation*}

%\medskip
%{\bf Moments of the two-direction multiscaling functions by doubling}

The main idea of the method by doubling in this section is to consider the
two-direction multiscaling function $\boldphi(x)$ and $\boldphi(-x)$ together,
and the two-direction multiwavelets $\boldpsi^{(s)}(x)$ and $\boldpsi^{(s)}(-x)$
together. That is, we consider refinable functions
$\boldPhi(x)=[\boldphi(x),\boldphi(-x)]^T$ and
$\boldPsi^{(s)}(x)=[\boldpsi^{(s)}(x),\boldpsi^{(s)}(-x)]^T$.
Then for $\boldPhi$ and $\boldPsi^{(s)}$, we develop a theory for computing
continuous moments. Finally, we choose parts of the continuous moments of
$\boldPhi$ and $\boldPsi^{(s)}$ for the continuous moments of $\boldphi$ and
$\boldpsi^{(s)}$, respectively for $s=1,2,\ldots,d-1$.

The $j$th {\em discrete moment} $M_{j}$ of the two-direction multiscaling function
$\boldphi$ is an $r\times r$ matrix defined by
\begin{equation}\label{dmphi0}
M_{j} = \frac{1}{\sqrt{d}} \sum_{k\in\Z} k^{j} [ P^{+}_k + P^{-}_k ], \qquad j=0,1,2,\ldots.
\end{equation}

The $j$th {\em discrete moment} $M^{\pm}_{j}$ of the one-direction multiscaling
function $\boldPhi$ is a $2r\times 2r$ matrix defined by
\begin{equation}\label{discmoment}
M^{\pm}_{j} = \frac{1}{\sqrt{d}} \sum_{k\in\Z} k^{j}
\begin{bmatrix}
    P^{+}_k & P^{-}_k \\ \noalign{\smallskip}
    P^{-}_{-k} & P^{+}_{-k} \\
\end{bmatrix}
, \qquad j=0,1,2,\ldots.
\end{equation}
%Discrete moments are related to the symbol by
%\begin{equation}\label{}
%    M^{\pm}_j = i^j D^j H(0),
%\end{equation}
%where $D^j$ stands for the $j$th derivative.
In particular,
\begin{equation*}
M^{\pm}_{0} = \frac{1}{\sqrt{d}} \sum_{k\in\Z}
\begin{bmatrix}
    P^{+}_k & P^{-}_k \\ \noalign{\smallskip}
    P^{-}_{-k} & P^{+}_{-k} \\
\end{bmatrix}
= \bold{P}_{\boldPhi}(1)
\end{equation*}
and $M^{\pm}_0$ is used in defining Condition E for $\boldphi$ in ~\eqref{condsymbol}.

%The symbol and the discrete moments are uniquely defined and easy to
%calculate. They are related by
%\begin{equation}
%%M^{\pm}{j} = i^{j} H^{(j)}(0),
%M^{\pm}_{j} = (-i)^{j} D^j H_{\boldPhi}(0),
%\end{equation}
%%where the superscript $(j)$ denotes the $j$th derivative.
%where the $D^j$ denotes the $j$th derivative.

The $j$th {\em continuous moment} $\boldm_{j}$ of $\boldphi$ is a vector of size $r$
defined by
\begin{equation}\label{cm}
    \boldm_{j} = \int_{-\infty}^{\infty} x^{j} \boldphi(x)\;\textrm{d}x, \qquad j=0,1,2,\ldots.
\end{equation}

The $j$th {\em continuous moment} $\boldm^{\pm}_{j}$ of $\boldPhi$ is a vector of
size $2r$ defined by
\begin{equation}\label{cmdef}
    \boldm^{\pm}_{j} = \int_{-\infty}^{\infty} x^{j} \boldPhi(x)\;\textrm{d}x, \qquad j=0,1,2,\ldots.
\end{equation}

In multiwavelet theory, computation of continuous moments for two-direction
multiwavelets %associated with two-direction multiscaling functions
is important, since vanishing continuous moments for two-direction multiwavelets
%associated with two-direction multiscaling functions
provide the approximation order for two-direction multiscaling functions.

Let $\boldpsi^{(s)}$, $s=1,2,\ldots,d-1$, be orthogonal two-direction multiwavelets
associated with the orthogonal two-direction multiscaling function $\boldphi$.
Let
\begin{equation} \label{defPsi}
\boldPsi^{(s)}(x)=[\boldpsi^{(s)}(x), \boldpsi^{(s)}(-x)]^T
\end{equation}
for $s=1,2,\ldots,d-1$.

The $j$th {\em discrete moments} $N^{(s)\pm}_{j}$ of $\boldPsi^{(s)}$, $s=1,2,\ldots,d-1$,
are $2r\times 2r$ matrices defined by
\begin{equation}\label{discmomentPsi}
N^{(s)\pm}_{j} = \frac{1}{\sqrt{d}} \sum_{k\in\Z} k^{j}
\begin{bmatrix}
    Q^{(s)+}_k & Q^{(s)-}_k \\ \noalign{\smallskip}
    Q^{(s)-}_{-k} & Q^{(s)+}_{-k} \\
\end{bmatrix}, \qquad j=0,1,2,\ldots.
\end{equation}

The $j$th {\em continuous moments} $\boldn^{(s)}_{j}$ of $\boldpsi^{(s)}$,
$s=1,2,\ldots,d-1$, are vectors of size $r$ defined by
\begin{equation}\label{cmpsi}
    \boldn^{(s)}_{j} = \int_{-\infty}^{\infty} x^{j} \boldpsi^{(s)}(x)\;\textrm{d}x,
    \qquad j=0,1,2,\ldots.
\end{equation}

The $j$th {\em continuous moments} $\boldn^{(s)\pm}_{j}$ of $\boldPsi^{(s)}$,
$s=1,2,\ldots,d-1$, are vectors of size $2r$ defined by
\begin{equation}\label{cmpsidef}
    \boldn^{(s)\pm}_{j} = \int_{-\infty}^{\infty} x^{j} \boldPsi^{(s)}(x)\;\textrm{d}x,
    \qquad j=0,1,2,\ldots.
\end{equation}

%For orthogonal two-direction multiwavelets $\boldpsi^{(s)}$ associated with
%orthogonal two-direction multiscaling function $\boldphi$,
%one can similarly compute discrete and continuous moments of $\boldpsi^{(s)}$
%in terms of the recurrence coefficients $Q^{(s)+}_k$, $Q^{(s)-}_k$ of
%$\boldpsi^{(s)}$ for $s=1,2,\ldots,d-1$.

The $j$th continuous moments $\boldm^{\pm}_{j}$ of $\boldPhi$ and
$\boldn^{(s)\pm}_{j}$ of $\boldPsi^{(s)}$ satisfy the following theorem.

\begin{theorem}\label{thm:PhiPsi}
Let the two-direction multiscaling functions $\boldphi$ and multiwavelets
$\boldpsi^{(s)}$ be defined as in ~\eqref{recrel} and ~\eqref{recrelpsi}, respectively,
and the refinable functions $\boldPhi$ and $\boldPsi^{(s)}$ defined as in
~\eqref{recrelPhi} and ~\eqref{defPsi}, respectively for $s=1,2,\ldots,d-1$.
Then the $j$th continuous moments $\boldm^{\pm}_{j}$ of $\boldPhi$ and
$\boldn^{(s)\pm}_{j}$ of $\boldPsi^{(s)}$ associated with $\boldPhi$ satisfy the
following:
\begin{equation}\label{mandM}
\boldm^{\pm}_{j} = d^{-j} \sum_{\ell=0}^{j} \binomial{j}{\ell} M^{\pm}_{j-\ell} \boldm^{\pm}_{\ell}, \qquad j=0,1,2,\ldots,
\end{equation}
\begin{equation}\label{nandN}
\boldn^{(s)\pm}_{j} = d^{-j} \sum_{\ell=0}^{j} \binomial{j}{\ell} N^{(s)\pm}_{j-\ell} \boldm^{\pm}_{\ell}, \qquad j=0,1,2,\ldots,
\end{equation}
where $\displaystyle \binomial{j}{\ell} = \dfrac{j!}{\ell!\,(j-\ell)!}$ stands
for the binomial coefficient.

%In particular,
%\begin{equation*}
%    \boldm^{\pm}_{0} = M^{\pm}_{0} \boldm^{\pm}_{0}.
%\end{equation*}
%
%Once $\boldm^{\pm}_{0}$ has been chosen, all other continuous moments are
%uniquely defined and can be computed recursively from these relations.
\end{theorem}

\begin{proof}
By substituting the recursion formula~\eqref{recrelPhi} into the integral in
~\eqref{cmdef}, we find after simplification that the continuous moments
$\boldm^{\pm}_{j}$ satisfy ~\eqref{mandM}.

By substituting the recursion formula~\eqref{recrelpsi} and %~\eqref{recrelneg}
$\boldpsi^{(s)}(-x)$ into the integral in ~\eqref{cmpsidef}, we find after
simplification that the continuous moments $\boldn^{(s)\pm}_{j}$ satisfy ~\eqref{nandN}.
\end{proof}

\begin{remark}\label{rmk:sphipsi}{\rm
We give some further explanatory remarks.

For $s=1,2,\ldots,d-1$ and $j=0,1,2,\ldots$,

1. By substitution, we obtain other formulas for $\boldm^{\pm}_j$ and $\boldn^{(s)\pm}_j$:
\begin{equation}\label{eq:cmPhiPsi2}
\boldm^{\pm}_{j} = \begin{bmatrix}
    1 \\ \noalign{\smallskip} (-1)^{j} \\
    \end{bmatrix}
    \boldm_j, \qquad
\boldn^{(s)\pm}_{j} = \begin{bmatrix}
    1 \\ \noalign{\smallskip} (-1)^{j} \\
    \end{bmatrix}
    \boldn^{(s)}_j.
\end{equation}

2. ~\eqref{eq:cmPhiPsi2} is useful for double checking the results of
$\boldm^{\pm}_j$ and $\boldn^{(s)\pm}_j$. % for $s=1,2,\ldots,d-1$.

3. Theorem \ref{thm:PhiPsi} provides one way to compute the continuous moments
$\boldm^{\pm}_j$ of $\boldPhi$ and $\boldn^{(s)\pm}_j$ of $\boldPsi^{(s)}$.
%for $s=1,2,\ldots,d-1$.
We will pursue this in the following.

%4. By~\eqref{eq:cmPhiPsi2},
%the upper half of $\boldm^{\pm}_{j}$ is $\boldm_j$,
%the lower half of $\boldm^{\pm}_{j}$ is $(-1)^{j} \boldm_j$,
%the upper half of $\boldn^{(s)\pm}_{j}$ is $\boldn^{(s)}_j$, and
%the lower half of $\boldn^{(s)\pm}_{j}$ is $(-1)^{j} \boldn^{(s)}_j$.
%%for $s=1,2,\ldots,d-1$.
}
\end{remark}

The normalizing condition for one-direction multiscaling functions $\boldphi$
is given in ~\cite{Keinert2004,Kwon2009}. In the following Lemma we derive
the normalizing condition for the two-direction multiscaling functions $\boldphi$.

\begin{theorem}\label{thm:phinormal} {\rm (Normalization for $\boldm_0$)}
Let $\boldphi$ be a two-direction orthogonal multiscaling function. Then the
normalizing condition for $\boldm_0$ of $\boldphi$ is
\begin{equation} \label{normal_m0}
\boldm_0^* \boldm_0 = \frac{1}{2}.
\end{equation}
\end{theorem}

\begin{proof}
By expanding the constant 1 in an orthogonal two-direction multiscaling
function series, we have
\begin{equation}\label{oneexpansion}
\begin{split}
1 &= \sum_{k\in\Z} [ \< 1, \boldphi(x-k) \> \,\boldphi(x-k)
   + \< 1, \boldphi(k-x) \> \,\boldphi(k-x) ] \\
&= \boldm_{0}^{*} \sum_{k\in\Z} [ \boldphi(x-k) + \boldphi(k-x) ].
\end{split}
\end{equation}
By integrating~\eqref{oneexpansion} on $[0,1]$, we have
\begin{equation*}
1 = \int_{0}^{1} 1 \;\textrm{d}x
= \boldm_{0}^{*} \sum_{k\in\Z} \int_{0}^{1} [ \boldphi(x-k) + \boldphi(k-x) ] \;\textrm{d}x
= \boldm_{0}^{*} \int_{-\infty}^{\infty} [ \boldphi(x) + \boldphi(-x) ] \;\textrm{d}x
= 2\boldm_{0}^{*} \boldm_{0}.
\end{equation*}
Hence, we have equation~\eqref{normal_m0}.
\end{proof}

%In orthogonal case the normalization $\boldm_0^* \boldm_0 = \frac{1}{2}$
%is mandatory.
For the scalar case $r=1$, that is, for orthogonal two-direction
scaling functions $\phi$, the normalization is
\begin{equation}\label{scalarmo}
    m_0 = \frac{\sqrt{2}}{2}.
\end{equation}

\begin{theorem}\label{thm:Phinormal}  {\rm (Normalization for $\boldm^{\pm}_0$)}
Let $\boldphi$ be a two-direction orthogonal multiscaling function.
Then the normalizing condition for $\boldm^{\pm}_0$ of $\boldPhi$ is
\begin{equation} \label{normal_mpm0}
(\boldm^{\pm}_0)^* \boldm^{\pm}_0 = 1.
\end{equation}
\end{theorem}

\begin{proof}
By ~\eqref{eq:cmPhiPsi2} and theorem ~\ref{thm:phinormal}, we have
\begin{equation*}
    (\boldm^{\pm}_0)^* \boldm^{\pm}_0
    = \boldm_0^* \begin{bmatrix}
                     1 & 1 \\
                   \end{bmatrix}
            \begin{bmatrix}
                     1  \\ 1 \\
                   \end{bmatrix}
            \boldm_0
    = 2 \boldm_0^* \boldm_0 = 1.
\end{equation*}
\end{proof}

From equation ~\eqref{recrelPhi}, $\boldPhi$ can be considered as a one-direction
multiscaling function. The result in theorem~\ref{thm:Phinormal} coincides with
the normalizing condition for the one-direction multiscaling functions which
is known as $(\boldm_0)^* \boldm_0 = 1$
%$(\boldm^{\pm}_0)^* \boldm^{\pm}_0 = 1$
(see~\cite{Keinert2004,Kwon2009}).

%Again in orthogonal case the normalization
%$(\boldm^{\pm}_0)^* \boldm^{\pm}_0 = 1$ is mandatory.

Now we summarize the above results in the following theorem.
The following theorem is the main result of this section.

\begin{theorem}\label{thm:sum}
Assume that a two-direction multiscaling function $\boldphi$ is orthogonal,
has compact support, is continuous (which implies approximation order
at least 1), and satisfies Condition E.
Let $\boldPhi$ be defined as in ~\eqref{recrelPhi}.
Then the continuous moments $\boldm^{\pm}_j$ of $\boldPhi$ and $\boldn^{(s)\pm}_j$ of $\boldPsi^{(s)}$ can be computed recursively as:\\
for $j=0$,
\begin{equation*}
\boldm^{\pm}_{0} = M^{\pm}_{0} \boldm^{\pm}_{0}, \qquad \text{with normalization}
\qquad (\boldm^{\pm}_0)^* \boldm^{\pm}_0 = 1;
\end{equation*}
for $j=1,2,3,\ldots$,
\begin{equation}\label{mcomputation2}
\boldm^{\pm}_{j} = \left( d^{j}I_{2r} - M^{\pm}_{0} \right)^{-1} \sum_{\ell=0}^{j-1}
\binomial{j}{\ell} M^{\pm}_{j-\ell} \boldm^{\pm}_{\ell};
\end{equation}
for $j=0,1,2,\ldots$,
\begin{equation}\label{nandN02}
\boldn^{(s)\pm}_{j} = d^{-j} \sum_{\ell=0}^{j} \binomial{j}{\ell} N^{(s)\pm}_{j-\ell} \boldm^{\pm}_{\ell},
\end{equation}
where $I_{2r}$ is the $2r \times 2r$ identity matrix.

Once $\boldm^{\pm}_{0}$ has been chosen, all other continuous moments are uniquely
defined and can be computed recursively from ~\eqref{mcomputation2} and ~\eqref{nandN02}.

Finally, the $j$th continuous moment $\boldm_{j}$ of $\boldphi$ and $\boldn^{(s)}_{j}$
of $\boldpsi^{(s)}$ are the upper half of $\boldm^{\pm}_{j}$ and $\boldn^{(s)\pm}_{j}$, respectively for $s=1,2,\ldots,d-1$.
\end{theorem}

\begin{proof}
By setting $j=0$ in equation ~\eqref{mandM}, we have
\begin{equation*}
    \boldm^{\pm}_{0} = M^{\pm}_{0} \boldm^{\pm}_{0}.
\end{equation*}
By theorem~\ref{thm:Phinormal}, we normalize $\boldm^{\pm}_{0}$ so that
$(\boldm^{\pm}_0)^* \boldm^{\pm}_0 = 1$.

From equation~\eqref{mandM}, we collect $\boldm^{\pm}_j$ terms on the left-hand side
and multiply by $d^j$ on both sides. Then the coefficient matrix
\begin{equation*}
    ( d^{j} I_{2r} - M^{\pm}_{0} )
\end{equation*}
of $\boldm^{\pm}_j$ is nonsingular by Condition E for $\boldphi$.
Hence, $\boldm^{\pm}_j$ can be determined recursively as
\begin{equation*}%\label{mandM22}
\boldm^{\pm}_{j} = \left( d^{j} I_{2r} - M^{\pm}_{0} \right)^{-1}
\sum_{\ell=0}^{j-1} \binomial{j}{\ell} M^{\pm}_{j-\ell} \boldm^{\pm}_{\ell}, \qquad j=1,2,3,\ldots.
\end{equation*}

\eqref{nandN02} is derived in ~\eqref{nandN}.
\end{proof}

%\begin{remark}{\rm
For the scalar case $r=1$, that is, for orthogonal two-direction scaling
function $\phi$, the zeroth continuous moment of $\boldPhi$ is
$$
\boldm^{\pm}_{0} = \frac{\sqrt{2}}{2} \begin{bmatrix}
    1 \\  1 \\
    \end{bmatrix}
$$
by remark ~\ref{rmk:sphipsi} and theorem ~\ref{thm:Phinormal}.
Hence, $m_{0} = \sqrt{2}/2$, which coincides with ~\eqref{scalarmo}.
%}
%\end{remark}

%\medskip
%{\bf Moments of two-direction multiwavelet by doubling}\\

%We will pursue this in the following.

%{\bf Algorithm for computing continuous moments $\boldm_j$ and $\boldn^{(s)}_j$ of
%$\boldphi$ and $\boldpsi^{(s)}$, respectively, for $j=0,1,2,\ldots$ by doubling.}
Now we provide an algorithm for computing continuous moments $\boldm_j$ and
$\boldn^{(s)}_j$ of $\boldphi$ and $\boldpsi^{(s)}$, respectively, for
$s=1,2,\ldots,d-1$ and $j=0,1,2,\ldots$ by doubling.

\begin{algorithm}\label{algo:mPhiPsi}{\rm
Let an orthogonal two-direction multiscaling function $\boldphi$ be given.
Let $\boldpsi^{(s)}$, $s=1,2,\ldots,d-1$, be orthogonal multiwavelets associated
with $\boldphi$.
Then an algorithm for computing continuous moments $\boldm_j$ of $\boldphi$
and $\boldn^{(s)}_j$ of $\boldpsi^{(s)}$ by doubling is the following:
\begin{itemize}
\item step 1: Compute the $j$th discrete moment $M^{\pm}_j$ of $\boldPhi$ for $j=0,1,2,\ldots$;
\item step 2: Compute the $j$th discrete moment $N^{(s)\pm}_j$ of $\boldPsi^{(s)}$ for $s=1,2,\ldots,d-1$ and $j=0,1,2,\ldots$;
\item step 3: Compute the zeroth continuous moment $\boldm^{\pm}_0$ of $\boldPhi$
    as the corresponding eigenvector to the eigenvalue 1 of $M^{\pm}_0$;
\item step 4: Normalize $\boldm^{\pm}_0$ so that $(\boldm^{\pm}_0)^* \boldm^{\pm}_0 = 1$;
\item step 5: Compute the $j$th continuous moment $\boldm^{\pm}_j$
    recursively for $j=1,2,\ldots$; %from equation ~\eqref{mandM}
\item step 6: Take the upper half of $\boldm^{\pm}_j$ as the $j$th continuous moments $\boldm_j$ of $\boldphi$ for $j=0,1,2,\ldots$;
\item step 7: Compute the $j$th continuous moment $\boldn^{(s)\pm}_j$
    recursively for $s=1,2,\ldots,d-1$ and $j=0,1,2,\ldots$; %from equation ~\eqref{nandN}
\item step 8: Take the upper half of $\boldn^{(s)\pm}_j$ as the $j$th continuous moment $\boldn^{(s)}_j$ of $\boldpsi^{(s)}$ for $s=1,2,\ldots,d-1$ and $j=0,1,2,\ldots$.
\end{itemize}
Step 2 can be placed anywhere between step 1 and step 7.
}
\end{algorithm}

The $j$th discrete moments $M_j$ of $\boldphi$ can be computed directly from
~\eqref{dmphi0}.

%We note that the method for computing moments by doubling is only good for
%the continuous moments, but not for the discrete moments, of $\boldphi$ and
%$\boldpsi^{(s)}$ for $s=1,2,\ldots,d-1$, since the discrete moments of $\boldphi$
%and $\boldpsi^{(s)}$ can not be computed by choosing parts of the discrete moments
%of $\boldPhi$ and $\boldPsi^{(s)}$, respectively for $s=1,2,\ldots,d-1$.
%This is not a problem in practice, because the discrete moments of $\boldphi$
%and $\boldpsi^{(s)}$ for $s=1,2,\ldots,d-1$ are just used in computing the continuous
%moments of $\boldphi$ and $\boldpsi^{(s)}$ for $s=1,2,\ldots,d-1$.
%(This will be more clear in next section.)
%Moreover, the discrete moments are not as important as the continuous moments
%in practice.

%%%%%%%%%%%%%%%%%%%%%%%%%%%%%%%%%%%%%%%%%%%%%%%%%%%%%%%%%%%%%%%%%%%%%%%%%%%%%%%%%%%
\bigskip %\medskip
\section{Moments of two-direction multiwavelets by separation} \label{sec:moments2}
%%%%%%%%%%%%%%%%%%%%%%%%%%%%%%%%%%%%%%%%%%%%%%%%%%%%%%%%%%%%%%%%%%%%%%%%%%%%%%%%%%%

In this section we provide another simple and efficient method for computing
the discrete and continuous moments of the two-direction multiscaling function and
two-direction multiwavelet associated with the two-direction multiscaling function.

The main idea of the method in this section is to separate the two-direction
multiscaling function $\boldphi$ and two-direction multiwavelet $\boldpsi$
%associated with $\boldphi$
into two parts, positive and negative, and then use the direct method on each part.

We recall ~\eqref{dmphi0} that the $j$th {\em discrete moment} $M_{j}$ of the two-direction multiscaling
function $\boldphi$ is an $r\times r$ matrix defined by
\begin{equation}\label{dmphi}
M_{j} = \frac{1}{\sqrt{d}} \sum_{k\in\Z} k^{j} [ P^{+}_k + P^{-}_k ], \qquad j=0,1,2,\ldots.
\end{equation}

Now we separate the $j$th discrete moment $M_j$ into positive and negative parts.
We define the $j$th {\em positive} and {\em negative discrete moments},
$M^{+}_{j}$ and $M^{-}_{j}$ respectively, of the discrete moment $M_{j}$ of
$\boldphi$ by
\begin{gather}\label{dmp2def}
    M^{+}_{j} = \frac{1}{\sqrt{d}} \sum_k k^j P^{+}_k, \qquad j=0,1,2,\ldots, \\
    M^{-}_{j} = \frac{1}{\sqrt{d}} \sum_k k^j P^{-}_k, \qquad j=0,1,2,\ldots,
\end{gather}
so that $M_{j} = M^{+}_{j} + M^{-}_{j}$.

By applying ~\eqref{recrel} in ~\eqref{cm}, the $j$th continuous moment $\boldm_{j}$
of two-direction multiscaling function $\boldphi$ can be expressed as
\begin{equation*}\label{cm2}
    \boldm_{j} = \sqrt{d} \sum_k \left[ P^{+}_k \int_{-\infty}^{\infty} x^{j} \boldphi(dx-k) \;\textrm{d}x + P^{-}_k \int_{-\infty}^{\infty} x^{j} \boldphi(k-dx)  \;\textrm{d}x \right], \qquad j=0,1,2,\ldots.
\end{equation*}

Now we separate the $j$th continuous moment $\boldm_{j}$ into positive and negative parts.
We define the $j$th {\em positive} and {\em negative continuous moments},
$\boldm^{+}_{j}$ and $\boldm^{-}_{j}$ respectively, of the continuous moment
$\boldm_{j}$ of $\boldphi$ by
\begin{gather}\label{cmp2def}
    \boldm^{+}_{j} =  \sqrt{d} \sum_k P^{+}_k \int_{-\infty}^{\infty} x^{j} \boldphi(dx-k)\;\textrm{d}x, \qquad j=0,1,2,\ldots, \\
    \boldm^{-}_{j} = \sqrt{d} \sum_k P^{-}_k \int_{-\infty}^{\infty} x^{j} \boldphi(k-dx)\;\textrm{d}x, \qquad j=0,1,2,\ldots,
\end{gather}
so that $\boldm_{j} = \boldm^{+}_{j} + \boldm^{-}_{j}$.

Let us define the $j$th {\em discrete moment} of the two-direction multiwavelet
$\boldpsi^{(s)}$, $s=1,2,\ldots,d-1$ by an $r\times r$ matrix
\begin{equation}\label{dmpsi}
N^{(s)}_{j} = \frac{1}{\sqrt{d}} \sum_{k\in\Z} k^{j} [ Q^{(s)+}_k + Q^{(s)-}_k ], \qquad j=0,1,2,\ldots.
\end{equation}

Now we separate the $j$th discrete moment $N_j$ into positive and negative parts.
We define the $j$th {\em positive} and {\em negative discrete moments},
$N^{(s)+}_{j}$ and $N^{(s)-}_{j}$ respectively, of the discrete moment
$N^{(s)}_{j}$ of $\boldpsi^{(s)}$, $s=1,2,\ldots,d-1$, by
\begin{gather}\label{dmp2def2}
    N^{(s)+}_{j} = \frac{1}{\sqrt{d}} \sum_k k^j Q^{(s)+}_k, \qquad j=0,1,2,\ldots, \\
    N^{(s)-}_{j} = \frac{1}{\sqrt{d}} \sum_k k^j Q^{(s)-}_k, \qquad j=0,1,2,\ldots,
\end{gather}
so that $N^{(s)}_{j} = N^{(s)+}_{j} + N^{(s)-}_{j}$.

We recall ~\eqref{cmpsi}, the $j$th {\em continuous moment} $\boldn^{(s)}_{j}$ of
the two-direction multiwavelet $\boldpsi^{(s)}$, $s=1,2,\ldots,d-1$, is a vector
of size $r$ defined by
\begin{equation}\label{cmpsi2}
\boldn^{(s)}_{j} = \int_{-\infty}^{\infty} x^j \boldpsi^{(s)}(x)\;\textrm{d}x, \qquad j=0,1,2,\ldots.
\end{equation}

%\medskip
%{\bf Moments of two-direction multiscaling function by separation}%\\

\begin{theorem}\label{thm:phipsi2}
Let $\boldphi$ be a two-direction multiscaling function defined as ~\eqref{recrel}
and assume that the zeroth discrete moment $M_0$ satisfies Condition E.
Let $\boldphi$ be a two-direction multiscaling function defined as ~\eqref{recrel} and $\boldpsi$ be two-direction multiwavelet associated with $\boldphi$. Then the $j$th continuous moments $\boldm_{j}$ of $\boldphi$ and $\boldn^{(s)}_{j}$ of $\boldpsi^{(s)}$, $s=1,2,\ldots,d-1$, associated with $\boldphi$ satisfy the following:
%can be determined recursively as
\begin{equation}\label{mandM2}
\boldm_{j} = d^{-j} \sum_{\ell=0}^{j} \binomial{j}{\ell} \left[ M^{+}_{j-\ell} + (-1)^\ell M^{-}_{j-\ell}  \right] \boldm_{\ell}, \qquad j=0,1,2,\ldots,
\end{equation}
\begin{equation}\label{nandN2}
\boldn^{(s)}_{j} = d^{-j} \sum_{\ell=0}^{j} \binomial{j}{\ell} \left[ N^{(s)+}_{j-\ell} + (-1)^\ell N^{(s)-}_{j-\ell}  \right] \boldm_{\ell}, \qquad j=0,1,2,\ldots,
\end{equation}
where $\displaystyle \binomial{j}{\ell} = \dfrac{j!}{\ell!\,(j-\ell)!}$ stands for the
binomial coefficient.

%In particular,
%\begin{equation*}
%    \boldm_{0} = M_{0} \boldm_{0}.
%\end{equation*}
%
%Moreover,
%\begin{equation}\label{mandM22}
%\boldm_{j} = \left( d^{j} I_{r} - \left[ M^{+}_{0} + (-1)^j M^{-}_{0}  \right] %\right)^{-1}
%\sum_{\ell=0}^{j-1} \binomial{j}{\ell} \left[ M^{+}_{j-\ell} + (-1)^j M^{-}_{j-\ell}  %\right] \boldm_{\ell}, \qquad j=1,2,\ldots,
%\end{equation}
%where $I_{r}$ is the identity matrix of size $r\times r$.
%
%Once $\boldm_{0}$ has been chosen, all other continuous moments are
%uniquely defined and can be computed recursively from these relations.
\end{theorem}

\begin{proof}
After substituting the recursion formula~\eqref{recrel} into the integral in
~\eqref{cmdef}, we separate the $j$th continuous moment $\boldm_j$ into two parts,
positive and negative continuous moments $\boldm^{+}_{j}$ and $\boldm^{-}_{j}$.

By changing variables for $dx-k$ in $\boldm^{+}_{j}$, we have
\begin{eqnarray*}
\boldm^{+}_{j} &=& \sqrt{d} \sum_{k\in\Z} P^{+}_k \int_{-\infty}^{\infty} x^{j} \boldphi(dx-k)\;\textrm{d}x \\
&=& d^{-j} \, \frac{1}{\sqrt{d}} \sum_{k\in\Z} P^{+}_{k} \int_{-\infty}^{\infty} (x+k)^j\,\boldphi(x) \;\textrm{d}x.
\end{eqnarray*}
By using the binomial expansion, we have
\begin{eqnarray*}
\boldm^{+}_{j} &=& d^{-j} \, \sum_{\ell=0}^j {j\choose \ell} \frac{1}{\sqrt{d}} \sum_{k\in\Z} k^{j-\ell} P^{+}_{k} \int_{-\infty}^{\infty} x^\ell \boldphi(x) \;\textrm{d}x \\
&=& d^{-j} \, \sum_{\ell=0}^j {j\choose \ell} M^{+}_{j-\ell} \boldm_\ell.
\end{eqnarray*}

By changing variables for $k-dx$ in $\boldm^{-}_{j}$, we have
\begin{eqnarray*}
\boldm^{-}_{j} &=& \sqrt{d} \sum_{k\in\Z} P^{-}_k \int_{-\infty}^{\infty} x^{j} \boldphi(k-dx)\;\textrm{d}x \\
&=& d^{-j} \, \frac{1}{\sqrt{d}} \sum_{k\in\Z} P^{-}_{k} \int_{\infty}^{-\infty} (-x+k)^j\,\boldphi(x) (-1)\;\textrm{d}x \\
&=& d^{-j} \, \frac{1}{\sqrt{d}} \sum_{k\in\Z} P^{-}_{k} \int_{-\infty}^{\infty} (-x+k)^j\,\boldphi(x) \;\textrm{d}x.
\end{eqnarray*}
By using the binomial expansion, we have
\begin{eqnarray*}
\boldm^{-}_{j} &=& d^{-j} \, \sum_{\ell=0}^j {j\choose \ell} (-1)^\ell \frac{1}{\sqrt{d}} \sum_{k\in\Z} k^{j-\ell} P^{-}_{k} \int_{-\infty}^{\infty} x^\ell \boldphi(x) \;\textrm{d}x \\
&=& d^{-j} \, \sum_{\ell=0}^j {j\choose \ell} (-1)^\ell M^{-}_{j-\ell} \boldm_\ell.
\end{eqnarray*}
Hence,
\begin{equation*}
\boldm_{j} = \boldm^{+}_{j} + \boldm^{-}_{j}
= d^{-j} \sum_{\ell=0}^{j} \binomial{j}{\ell} \left[ M^{+}_{j-\ell} + (-1)^\ell M^{-}_{j-\ell}  \right] \boldm_{\ell}, \qquad j=0,1,2,\ldots.
\end{equation*}

%In particular,
%\begin{equation*}
%    \boldm_{0} = M_{0} \boldm_{0}.
%\end{equation*}
%
%For $j=1,2,\ldots$, let us collect $\boldm_j$ terms on the left-hand side
%and multiply by $d^j$ on both sides. Then the coefficient matrix
%\begin{equation*}
%    ( d^{j} I_{r} - \left[ M^{+}_{0} + (-1)^j M^{-}_{0} \right] )
%\end{equation*}
%of $\boldm_j$ is nonsingular by Condition E for $\boldphi$.
%Hence, $\boldm_j$ can be determined recursively as
%\begin{equation*}%\label{mandM22}
%\boldm_{j} = \left( d^{j} I_{r} - \left[ M^{+}_{0} + (-1)^j M^{-}_{0}  \right] %\right)^{-1}
%\sum_{\ell=0}^{j-1} \binomial{j}{\ell} \left[ M^{+}_{j-\ell} + (-1)^\ell M^{-}_{j-\ell}  %\right] \boldm_{\ell}, \qquad j=0,1,2,\ldots.
%\end{equation*}
%
%Once $\boldm_{0}$ has been chosen, all other continuous moments are
%uniquely defined and can be computed recursively from these relations.

For $\boldn^{(s)}_j$, it is similar to the proof of $\boldm_j$.
\end{proof}

We have proved the normalization
\begin{equation*}
    \boldm_{0}^* \,\boldm_{0} = \frac{1}{2}
\end{equation*}
in theorem ~\ref{thm:phinormal}.
Hence, it is not necessary to derive normalization here again.

Now we summarize the above results in the following theorem.
The following theorem is the main result of this section.

\begin{theorem}\label{thm:sum}
Assume that a two-direction multiscaling function $\boldphi$ is orthogonal,
has compact support, is continuous (which implies approximation order
at least 1), and satisfies Condition E.
Then the continuous moments $\boldm_j$ of $\boldphi$ and $\boldn^{(s)}_j$ of $\boldpsi^{(s)}$ can be computed recursively as:\\
for $j=0$,
\begin{equation*}
\boldm_{0} = M_{0} \,\boldm_{0}, \qquad \text{with normalization}
\qquad \boldm_0^* \,\boldm_0 = \frac{1}{2};
\end{equation*}
for $j=1,2,3,\ldots$,
\begin{equation}\label{mandM22}
\boldm_{j} = \left( d^{j}I_{r} - \left[ M^{+}_{0} + (-1)^j M^{-}_{0}  \right] \right)^{-1} \sum_{\ell=0}^{j-1}
\binomial{j}{\ell} M_{j-\ell} \boldm_{\ell};
\end{equation}
for $j=0,1,2,\ldots$,
\begin{equation}\label{nandN22}
\boldn^{(s)}_{j} = d^{-j} \sum_{\ell=0}^{j} \binomial{j}{\ell} \left[ N^{(s)+}_{j-\ell} + (-1)^\ell N^{(s)-}_{j-\ell}  \right] \boldm_{\ell},
\end{equation}
where $I_{r}$ is the $r \times r$ identity matrix.

Once $\boldm_{0}$ has been chosen, all other continuous moments are uniquely defined
and can be computed recursively from ~\eqref{mandM22} and ~\eqref{nandN22}.
%these relations.
\end{theorem}

\begin{proof}
By setting $j=0$ in equation ~\eqref{mandM2}, we have
\begin{equation*}
    \boldm_{0} = M_{0} \,\boldm_{0}.
\end{equation*}
By theorem~\ref{thm:phinormal}, we normalize $\boldm_{0}$ so that
$\boldm_0^* \,\boldm_0 = 1/2$.

From equation~\eqref{mandM2}, we collect $\boldm_j$ terms on the left-hand side
and multiply by $d^j$ on both sides. Then the coefficient matrix
\begin{equation*}
    ( d^{j} I_{r} - \left[ M^{+}_{0} + (-1)^j M^{-}_{0} \right] )
\end{equation*}
of $\boldm_j$ is nonsingular by Condition E for $\boldphi$.
Hence, $\boldm_j$ can be determined recursively as
\begin{equation*}%\label{mandM22}
\boldm_{j} = \left( d^{j} I_{r} - \left[ M^{+}_{0} + (-1)^j M^{-}_{0} \right] \right)^{-1}
\sum_{\ell=0}^{j-1} \binomial{j}{\ell} M_{j-\ell} \boldm_{\ell}, \qquad j=1,2,3,\ldots.
\end{equation*}

\eqref{nandN22} is derived in ~\eqref{nandN2}.
\end{proof}

For future reference, we list some in detail.
\begin{eqnarray*}
\boldm_{0} &=& M_0 \,\boldm_{0} \qquad \text{with normalization} \qquad \boldm_{0}^* \, \boldm_{0}=\frac{1}{2}, \\
\boldm_{1} &=& \left( d \,I_{r} - M^{+}_0 + M^{-}_0 \right)^{-1} M_1 \,\boldm_{0}, \\
\boldm_{2} &=& \left( d^{2} \,I_{r} - M_0 \right)^{-1} [ M_2 \,\boldm_{0} + 2(M^{+}_1 - M^{-}_1 )\,\boldm_{1} ], \\
\boldm_{3} &=& \left( d^{3} \,I_{r} - M^{+}_0 + M^{-}_0 \right)^{-1} [ M_3 \,\boldm_{0} + 3(M^{+}_2 - M^{-}_2 )\,\boldm_{1} + 3M_1\,\boldm_{2} ].
\end{eqnarray*}

%For future reference, we list some in detail.
For $s=1,2,\ldots,d-1$,
\begin{eqnarray*}
\boldn^{(s)}_{0} &=& N^{(s)}_0 \,\boldm_{0}, \\
\boldn^{(s)}_{1} &=& d^{-1} \,[ N^{(s)}_1 \,\boldm_{0} + (N^{(s)+}_0 - N^{(s)-}_0 )\,\boldm_{1} ], \\
\boldn^{(s)}_{2} &=& d^{-2} \,[ N^{(s)}_2 \,\boldm_{0} + 2(N^{(s)+}_1 - N^{(s)-}_1 )\,\boldm_{1} + N_0 \,\boldm_{2} ], \\
\boldn^{(s)}_{3} &=& d^{-3} \,[ N^{(s)}_3 \,\boldm_{0} + 3(N^{(s)+}_2 - N^{(s)-}_2 )\,\boldm_{1} + 3N^{(s)}_1\,\boldm_{2} + (N^{(s)+}_0 - N^{(s)-}_0 )\,\boldm_{3} ].
\end{eqnarray*}

%\medskip
%{\bf Moments of two-direction multiwavelet by separation}\\
%Now we work on moments for multiwavelets.

%{\bf Algorithm for computing continuous moments $\boldm_j$ and $\boldn^{(s)}_j$ of
%$\boldphi$ and $\boldpsi^{(s)}$, respectively, for $j=0,1,2,\ldots$ by separation.}
Now we provide an algorithm for computing continuous moments $\boldm_j$ and
$\boldn^{(s)}_j$ of $\boldphi$ and $\boldpsi^{(s)}$, respectively,
for $j=0,1,2,\ldots$ by separation.

\begin{algorithm}\label{algo:mphipsi2}{\rm
Let an orthogonal two-direction multiscaling function $\boldphi$ be given.
Let $\boldpsi^{(s)}$, $s=1,2,\ldots,d-1$, be orthogonal multiwavelets associated
with $\boldphi$.
Then an algorithm for computing continuous moments $\boldm_j$ of $\boldphi$ and $\boldn^{(s)}_j$ of $\boldpsi^{(s)}$ by separation is the following:
\begin{itemize}
\item step 1: Compute the $j$th positive and negative discrete moments, $M^{+}_j$ and $M^{-}_j$ respectively, and the $j$th discrete moment $M_j$ of $\boldphi$ for $j=0,1,2,\ldots$;
\item step 2: Compute the $j$th positive and negative discrete moments, $N^{(s)+}_j$ and $N^{(s)-}_j$ respectively, and the $j$th discrete moment $N^{(s)}_j$ of $\boldpsi^{(s)}$ for $s=1,2,\ldots,d-1$ and $j=0,1,2,\ldots$;
\item step 3: Compute the zeroth continuous moment $\boldm_0$ of $\boldphi$
    as the corresponding eigenvector to the eigenvalue 1 of $M_0$; %Eigenvector $\boldm_0$ of size $r$ has to be of the form $c[1,1,\ldots,1]^T$ for some constant $c$;
\item step 4: Normalize $\boldm_0$ so that $\boldm_0^* \boldm_0 = 1/2$;
\item step 5: Compute the $j$th continuous moment $\boldm_j$ of $\boldphi$ for $j=1,2,\ldots$;
\item step 6: Compute the $j$th continuous moments $\boldn^{(s)}_j$ of $\boldpsi^{(s)}$ for $s=1,2,\ldots,d-1$ and $j=0,1,2,\ldots$.
\end{itemize}
Step 2 can be placed anywhere between step 1 and step 6.
}
\end{algorithm}

\begin{remark}\label{}{\rm
We give some further explanatory remarks.

1. Two methods for computing continuous moments for $\boldphi$ and $\boldpsi^{(s)}$
for $s=1,2,\ldots,d-1$ by doubling and by separation generate the same results.
This will be demonstrated in Example~\ref{ex1}. %in section ~\ref{examples}.

2. The discrete and continuous moments for $\boldphi$ and $\boldpsi^{(s)}$
for $s=1,2,\ldots,d-1$ can be computed by both doubling and separation.

3. The method for computing moments by separation in section~\ref{sec:moments2} is
simpler, faster, and more efficient than that by doubling in section~\ref{moments}.
So, it is recommended to use the method for computing continuous moments
by separation in practice.

}
\end{remark}

%%%%%%%%%%%%%%%%%%%%%%%%%%%%%%%%%%%%%%%%%%%%%%%%%%%%%%%%%%%%%%%%%%%%%%%%%%%%%%%%%%%
\bigskip
\section{Examples}\label{examples}
%%%%%%%%%%%%%%%%%%%%%%%%%%%%%%%%%%%%%%%%%%%%%%%%%%%%%%%%%%%%%%%%%%%%%%%%%%%%%%%%%%%

In this section, we provide two examples for illustrating the general theory
in sections 2, 3 and 4.
%However, two methods for computing moments generate
%the same results, the method for computing moments by separation is simpler,
%faster, and more efficient than that by doubling.
%in the sense of computational effort and cost.

\bigskip
% Example 1
\begin{example} \label{ex1} {\rm
In this example we take two-direction scaling function $\phi$ and two-direction
wavelet $\psi$ associated with $\phi$ given in ~\cite{YangXie2008}.
We note that multiplicity $r=1$ and dilation factor $d=2$.

The nonzero recursion coefficients for $\phi$ are
\begin{gather*}
    P^{+}_{1} = \frac{1}{\sqrt{2}} \frac{3}{4}, \quad
    P^{+}_{2} = \frac{1}{\sqrt{2}} \frac{2-\sqrt{7}}{4}, \quad
    P^{-}_{2} = \frac{1}{\sqrt{2}} \frac{2+\sqrt{7}}{4}, \quad
    P^{-}_{3} = \frac{1}{\sqrt{2}} \frac{1}{4}.
\end{gather*}
The nonzero recursion coefficients for $\psi$ are
\begin{gather*}
    Q^{+}_{-2} = \frac{1}{\sqrt{2}} \frac{3}{4}, \quad
    Q^{+}_{-3} = -\frac{1}{\sqrt{2}} \frac{2-\sqrt{7}}{4}, \quad
    Q^{-}_{-1} = -\frac{1}{\sqrt{2}} \frac{2+\sqrt{7}}{4}, \quad
    Q^{-}_{-2} = \frac{1}{\sqrt{2}} \frac{1}{4}.
\end{gather*}
These differ from Yang and Xie ~\cite{YangXie2008} by a factor of
$1/\sqrt{2}$, due to differences in notation.

$\phi$ is supported on $[0,2]$ and $\psi$ associated with $\phi$ is supported
on $[-2,0]$. $\phi$ is orthogonal and $\psi$ is also orthogonal.

The matrix
\begin{align*}\label{condsymbol}
\begin{bmatrix}
        P^{+}(1) & P^{-}(1) \\
        P^{-}(1) & P^{+}(1) \\
    \end{bmatrix}
= \frac{1}{\sqrt{d}} \sum_{k\in\Z}
    \begin{bmatrix}
        P^{+}_k & P^{-}_k \\ \noalign{\medskip}
        P^{-}_{k} & P^{+}_{k} \\
    \end{bmatrix}
= \frac{1}{8}
    \begin{bmatrix}
        5-\sqrt{7} & 3+\sqrt{7} \\ \noalign{\medskip}
        3+\sqrt{7} & 5-\sqrt{7} \\
    \end{bmatrix}
\end{align*}
has two eigenvalues 1 and $(1-\sqrt{7})/4 \approx -0.4114$. Hence, Condition E for $\phi$ is satisfied.

\medskip
{\bf Moments of two-direction multiwavelets by doubling.}\\
The discrete moments $M^{\pm}_j$ of $\boldPhi(x) = [\phi(x), \phi(-x)]^T$
for $j=0,1,2$ are
\begin{gather*}
M^{\pm}_{0} = \frac{1}{8}
        \begin{bmatrix}
        5-\sqrt{7} & 3+\sqrt{7} \\ \noalign{\smallskip}
        3+\sqrt{7} & 5-\sqrt{7}
        \end{bmatrix}, \quad
M^{\pm}_{1} = \frac{1}{8}
        \begin{bmatrix}
        7-2\sqrt{7} & 7+2\sqrt{7} \\ \noalign{\smallskip}
        -7-2\sqrt{7}  & -7+2\sqrt{7}
        \end{bmatrix}, \quad \\
M^{\pm}_{2} = \frac{1}{8}
        \begin{bmatrix}
        11-4\sqrt{7} & 17+4\sqrt{7}  \\ \noalign{\smallskip}
        17+4\sqrt{7}  & 11-4\sqrt{7}
        \end{bmatrix}.
\end{gather*}

The eigenvalues of $M^{\pm}_0$ are 1 and $(1-\sqrt{7})/4$, and the corresponding
eigenvector to the eigenvalue 1 is $[1, 1]^T$.
Normalization factor is $\sqrt{2}/2$.
Normalized $\boldm^{\pm}_0$ is
$\sqrt{2}/2 \begin{bmatrix} 1 & 1 \\ \end{bmatrix}^T$.
We find that $2^{j} I_{2} -  M^{\pm}_0$ is invertible for $j=1,2,3,\ldots$.
Hence, we can find the $j$th continuous moments $\boldm^{\pm}_j$ of $\Phi$
recursively, for $j=1,2,3,\ldots$.

The continuous moments $\boldm^{\pm}_j$ of $\boldPhi$ for $j=0,1,2$ are
\begin{gather*}
\boldm^{\pm}_{0} = \frac{\sqrt{2}}{2}
        \begin{bmatrix}
        1 \\ 1
        \end{bmatrix}, \quad
\boldm^{\pm}_{1} = \frac{7\sqrt{2}-\sqrt{14}}{12} \begin{bmatrix}
        1 \\ -1
        \end{bmatrix}, \quad
\boldm^{\pm}_{2} = \frac{28\sqrt{2}-7\sqrt{14}}{36} \begin{bmatrix}
        1 \\ 1
        \end{bmatrix}.
\end{gather*}
Hence, the $j$th continuous moments $m_j$ of $\phi$ for $j=0,1,2$ are
\begin{equation*}
m_0 = \frac{\sqrt{2}}{2},\qquad
m_1 = \frac{7\sqrt{2}-\sqrt{14}}{12}, \qquad
m_2 = \frac{28\sqrt{2}-7\sqrt{14}}{36}.
\end{equation*}

%We note that we have obtained $y^{+}_{0} = m_{0}$, $y^{-}_{0} = m_{0}$,
%$y^{+}_{1} = m_{1}$, and $y^{-}_{1} = - m_{1}$ as expected.

The discrete moments $N^{\pm}_j$ of $\boldPsi(x) = [\psi(x), \psi(-x)]^T$
for $j=0,1,2$ are
\begin{gather*}
N^{\pm}_{0} = \frac{1}{8}
        \begin{bmatrix}
        1+\sqrt{7} & -1-\sqrt{7} \\ \noalign{\smallskip}
        -1-\sqrt{7} & 1+\sqrt{7}
        \end{bmatrix}, \quad
N^{\pm}_{1} = \frac{\sqrt{7}}{8}
        \begin{bmatrix}
        -3 & 1 \\ \noalign{\smallskip}
        -1  & 3
        \end{bmatrix}, \quad %\\
N^{\pm}_{2} = \frac{1}{8}
        \begin{bmatrix}
        -6+9\sqrt{7} & 2-\sqrt{7}  \\ \noalign{\smallskip}
        2-\sqrt{7}  & -6+9\sqrt{7}
        \end{bmatrix}.
\end{gather*}

The continuous moments $\boldn^{\pm}_j$ of $\boldPsi$ for $j=0,1,2$ are
\begin{gather*}
\boldn^{\pm}_{0} =
        \begin{bmatrix}
        0 \\ 0
        \end{bmatrix}, \quad
\boldn^{\pm}_{1} = \begin{bmatrix}
        0 \\ 0
        \end{bmatrix}, \quad
\boldn^{\pm}_{2} =  \frac{4\sqrt{2}-\sqrt{14}}{48} \begin{bmatrix}
        1 \\ 1
        \end{bmatrix}.
\end{gather*}

Hence, the $j$th continuous moments $n_j$ of $\psi$ for $j=0,1,2$ are
\begin{gather*}
    n_0 = 0,\qquad
    n_1 = 0, \qquad
    n_2 = \frac{4\sqrt{2}-\sqrt{14}}{48}.
\end{gather*}

\medskip
{\bf Moments of two-direction multiwavelets by separation.}\\
Positive and negative discrete moments, $M^{+}_j$ and $M^{-}_j$ respectively,
of $\phi$ for $j=0,1,2$ are
\begin{equation*}
\begin{aligned}
    M^{+}_{0} &= \frac{5-\sqrt{7}}{8}, \\
    M^{-}_{0} &= \frac{3+\sqrt{7}}{8}, \\
\end{aligned}
\qquad
\begin{aligned}
    M^{+}_{1} &= \frac{7-2\sqrt{7}}{8}, \\
    M^{-}_{1} &= \frac{7+2\sqrt{7}}{8}, \\
\end{aligned}
\qquad
\begin{aligned}
    M^{+}_{2} &= \frac{11-4\sqrt{7}}{8}, \\
    M^{-}_{2} &= \frac{17+4\sqrt{7}}{8}.
\end{aligned}
\end{equation*}
Hence, the discrete moments $M_j$ of $\phi$ for $j=0,1,2$ are
\begin{gather*}
    M_{0} = 1, \quad
    M_{1} = \frac{7}{4}, \quad
    M_{2} = \frac{7}{2}.
\end{gather*}

Without computing the positive and negative continuous moments, $m^{+}_j$ and
$m^{-}_j$ respectively, of $\phi$, we can compute the $j$th continuous moments $m_j$ of $\phi$ for $j=0,1,2,\ldots$ by equation ~\eqref{mandM22}.

The fact that $M_0=1$ implies that the eigenvalue of $M_0$ is 1 and the
corresponding eigenvector to the eigenvalue 1 is 1. Normalization factor is $\sqrt{2}/2$. Hence, $m_0 = \sqrt{2}/2$.
We find that
\begin{equation*}
2^{j} - \left[ M^{+}_{0} + (-1)^j M^{-}_{0} \right]
= \begin{cases}
        2^{j}-\frac{1-\sqrt{7}}{4}, & \text{for odd $j$},\\ \noalign{\smallskip}
        2^{j}-1, \qquad & \text{for even $j$}.
    \end{cases}
\end{equation*}
$2^{j} - \left[ M^{+}_{0} + (-1)^j M^{-}_{0} \right]$
is not 0 and therefore invertible for $j=1,2,3,\ldots$.
Hence, we can find the $j$th continuous moments $m_j$ of $\phi$ recursively,
for $j=1,2,3,\ldots$.

The $j$th continuous moments $m_j$ of $\phi$ for $j=0,1,2$ are
\begin{equation*}
m_0 = \frac{\sqrt{2}}{2},\qquad
m_1 = \frac{7\sqrt{2}-\sqrt{14}}{12}, \qquad
m_2 = \frac{28\sqrt{2}-7\sqrt{14}}{36}.
\end{equation*}

The positive and negative discrete moments, $N^{+}_j$ and $N^{-}_j$ respectively,
of $\psi$ for $j=0,1,2$ are
\begin{equation*}
\begin{aligned}
    N^{+}_{0} &= \frac{1+\sqrt{7}}{8}, \\
    N^{-}_{0} &= -\frac{1+\sqrt{7}}{8}, \\
\end{aligned}
\qquad
\begin{aligned}
    N^{+}_{1} &= -\frac{3\sqrt{7}}{8}, \\
    N^{-}_{1} &= \frac{\sqrt{7}}{8}, \\
\end{aligned}
\qquad
\begin{aligned}
    N^{+}_{2} &= \frac{-6+9\sqrt{7}}{8}, \\
    N^{-}_{2} &= \frac{2-\sqrt{7}}{8}.
\end{aligned}
\end{equation*}
Hence, the discrete moments $N_j$ of $\psi$ for $j=0,1,2$ are
\begin{gather*}
    N_{0} = 0, \quad
    N_{1} = -\frac{\sqrt{7}}{4}, \quad
    N_{2} = -\frac{1}{2}+\sqrt{7}.
\end{gather*}

%Without computing the positive and negative continuous moments, $n^{+}_j$ and
%$n^{-}_j$ respectively, of $\psi$, we can compute the $j$th continuous
%moments $n_j$ of $\psi$ for $j=0,1,2,\ldots$ by equation ~\eqref{nandN22}.

The $j$th continuous moments $n_j$ of $\psi$ for $j=0,1,2$ are
\begin{gather*}
n_0 = 0, \qquad
n_1 = 0, \qquad
n_2 = \frac{4\sqrt{2}-\sqrt{14}}{48}.
\end{gather*}

Hence, two methods by doubling and by separation generate the same results
for the continuous moments $m_j$ and $n_j$ of $\phi$ and $\psi$,
respectively for $j=0,1,2,\ldots$.

Since $n_0 = 0 = n_1$ and $n_2 \neq 0$, $\psi$ has 2 vanishing
moments and $\phi$ provides approximation order 2.

}
$\hfill\Box$
\end{example} % end of example 1

% Example 2
\begin{example} \label{ex2} {\rm
In this example we take two-direction multiscaling function $\boldphi$ and two-direction
multiwavelet $\boldpsi$ associated with $\boldphi$ given in ~\cite{WangZhouWang2011}.
We note that multiplicity $r=2$ and dilation factor $d=2$.

The nonzero recursion coefficients for $\boldphi$ are
\begin{gather*}
    P^{+}_{-3} = \frac{1}{8\sqrt{2}} \begin{bmatrix}
                                     0 & 0 \\
                                     -2\sqrt{3}+\sqrt{21} & 0 \\
                                   \end{bmatrix}, \quad
    P^{+}_{-2} = \frac{1}{8\sqrt{2}} \begin{bmatrix}
                                     0 & 0 \\
                                    3\sqrt{3} & 0 \\
                                   \end{bmatrix}, \\
    P^{+}_{1} = \frac{1}{8\sqrt{2}} \begin{bmatrix}
                                     6 & 0 \\
                                     0 & 3 \\
                                   \end{bmatrix}, \quad
    P^{+}_{2} = \frac{1}{8\sqrt{2}} \begin{bmatrix}
                                     4-2\sqrt{7} & 0 \\
                                    0 & 2-\sqrt{7} \\
                                   \end{bmatrix},
\end{gather*}
and
\begin{gather*}
    P^{-}_{-2} = \frac{1}{8\sqrt{2}} \begin{bmatrix}
                                     0 & 0 \\
                                     \sqrt{3} & 0 \\
                                   \end{bmatrix}, \quad
    P^{-}_{-1} = \frac{1}{8\sqrt{2}} \begin{bmatrix}
                                     0 & 0 \\
                                    -2\sqrt{3}-\sqrt{21} & 0 \\
                                   \end{bmatrix}, \\
    P^{-}_{2} = \frac{1}{8\sqrt{2}} \begin{bmatrix}
                                     4+2\sqrt{7} & 0 \\
                                     0 & 2+\sqrt{7} \\
                                   \end{bmatrix}, \quad
    P^{-}_{3} = \frac{1}{8\sqrt{2}} \begin{bmatrix}
                                     2 & 0 \\
                                     0 & 1 \\
                                   \end{bmatrix}.
\end{gather*}
The nonzero recursion coefficients for $\boldpsi$ are
\begin{gather*}
    Q^{+}_{-3} = \frac{1}{8\sqrt{2}} \begin{bmatrix}
                                     0 & -4+2\sqrt{7} \\
                                     -2+\sqrt{7} & 0 \\
                                   \end{bmatrix}, \quad
    Q^{+}_{-2} = \frac{1}{8\sqrt{2}} \begin{bmatrix}
                                     0 & 6 \\
                                     3 & 0 \\
                                   \end{bmatrix}, \\
    Q^{+}_{1} = \frac{1}{8\sqrt{2}} \begin{bmatrix}
                                     0 & 0 \\
                                     0 & -3\sqrt{3} \\
                                   \end{bmatrix}, \quad
    Q^{+}_{2} = \frac{1}{8\sqrt{2}} \begin{bmatrix}
                                     0 & 0 \\
                                    0 & -2\sqrt{3}+\sqrt{21} \\
                                   \end{bmatrix},
\end{gather*}
and
\begin{gather*}
    Q^{-}_{-2} = \frac{1}{8\sqrt{2}} \begin{bmatrix}
                                     0 & 2 \\
                                     1 & 0 \\
                                   \end{bmatrix}, \quad
    Q^{-}_{-1} = \frac{1}{8\sqrt{2}} \begin{bmatrix}
                                     0 & -4-2\sqrt{7} \\
                                    -2-\sqrt{7} & 0 \\
                                   \end{bmatrix}, \\
    Q^{-}_{2} = \frac{1}{8\sqrt{2}} \begin{bmatrix}
                                     0 & 0 \\
                                     0 & -2\sqrt{3}-\sqrt{21} \\
                                   \end{bmatrix}, \quad
    Q^{-}_{3} = \frac{1}{8\sqrt{2}} \begin{bmatrix}
                                     0 & 0 \\
                                     0 & -\sqrt{3} \\
                                   \end{bmatrix}.
\end{gather*}
These differ from Wang, Zhou, and Wang ~\cite{WangZhouWang2011} by a factor
of $1/\sqrt{2}$, due to differences in notation.

$\phi_1$ is supported on $[0,2]$ and support of $\phi_2$ is contained in $[-3,3]$.
%$\boldpsi$ associated with $\boldphi$ is supported
$\psi_1$ is supported on $[-2,0]$ and support of $\psi_2$ is contained in $[-3,3]$.

$\boldphi$ is orthogonal and $\boldpsi$ is also orthogonal.

The matrix
\begin{align*}\label{condsymbol}
\begin{bmatrix}
        P^{+}(1) & P^{-}(1) \\
        P^{-}(1) & P^{+}(1) \\
    \end{bmatrix}
= \frac{1}{\sqrt{d}} \sum_{k\in\Z}
    \begin{bmatrix}
        P^{+}_k & P^{-}_k \\ \noalign{\medskip}
        P^{-}_{k} & P^{+}_{k} \\
    \end{bmatrix}
\end{align*}
has four eigenvalues $1$, $1/2$, $(3-3\sqrt{7}+\sqrt{8-2\sqrt{7}})/16 \approx -0.2057$
and $(3-3\sqrt{7}-\sqrt{8-2\sqrt{7}})/16 \approx -0.4114$. \
Hence, Condition E for $\boldphi$ is satisfied.

%$\boldPhi$ provides approximation order 2 with approximation vectors
%\begin{equation*}
%   \frac{\sqrt{2}}{2} \begin{bmatrix}
%     1  \\ 0 \\
%    \end{bmatrix}, \qquad
%   \frac{7\sqrt{2}-\sqrt{14}}{12} \begin{bmatrix}
%     1  \\ -1 \\
%    \end{bmatrix}.
%\end{equation*}
%Hence, $\boldphi$ provides approximation order 2 with approximation vectors
%\begin{equation*}
%   y^{+}_{0} = \frac{\sqrt{2}}{2}, \qquad y^{-}_{0} = \frac{\sqrt{2}}{2}, \qquad
%   y^{+}_{1} = \frac{7\sqrt{2}-\sqrt{14}}{12}, \qquad y^{-}_{1} = - %\frac{7\sqrt{2}-\sqrt{14}}{12}.
%\end{equation*}

%% Graph was created by running ex2pointvalue.m
%% Graph was created by running multipointvalue.m
%\begin{figure}[h] \label{fig:wzwPhiPsi}  \centering
%\includegraphics[width=6.6in,height=2.6in]{graphics/wzwPhiPsi.jpg}
%\caption{
%Graph (left) of orthogonal two-direction multiscaling function $\boldphi=[\phi, %\frac{\sqrt{3}}{2}\psi + \frac{1}{2}\phi]^T$, where $\phi$ and $\psi$ are in
%example 1. $\boldphi$ is supported on $[-2,2]$.  $\boldphi$ provides approximation order %2.
%Graph (right) of orthogonal multiwavelet $\boldpsi$ associated with $\boldpsi=[\psi, %\frac{1}{2}\psi - \frac{\sqrt{3}}{2}\phi]^T$, where $\phi$ and $\psi$ are in
%example 1. $\boldpsi$ is supported on $[-2,2]$. $\boldpsi$ has 2 vanishing moments.
%}
%\end{figure}

In this example, we only compute moments by separation for simplicity.

\medskip
{\bf Moments of two-direction multiwavelets by separation.}\\
The positive and negative zeroth discrete moments, $M^{+}_0$ and $M^{-}_0$, of
$\boldphi$ are
\begin{gather*}
    M^{+}_{0} = \frac{1}{16} \begin{bmatrix}
                   10-2\sqrt{7} & 0 \\
                   \sqrt{3}+\sqrt{21} & 5-\sqrt{7} \\
                 \end{bmatrix}, \quad
    M^{-}_{0} = \frac{1}{16} \begin{bmatrix}
                   2\sqrt{7}+6 & 0 \\
                   -\sqrt{3}-\sqrt{21} & 3+\sqrt{7} \\
                 \end{bmatrix}.
\end{gather*}
The discrete moments, $M_j$ of $\boldphi$ for $j=0,1,2$ are
\begin{gather*}
    M_{0} = \frac{1}{2} \begin{bmatrix}
                   2 & 0 \\
                   0 & 1 \\
                 \end{bmatrix}, \quad
    M_{1} = \frac{1}{8} \begin{bmatrix}
                   14 & 0 \\
                   -\sqrt{21} & 7 \\
                 \end{bmatrix}, \quad
    M_{2} = \frac{1}{4} \begin{bmatrix}
                   14 & 0 \\
                   -\sqrt{3}+2\sqrt{21} & 7 \\
                 \end{bmatrix}.
\end{gather*}

%Without computing the positive and negative continuous moments, $\boldm^{+}_j$ and
%$\boldm^{-}_j$ respectively, of $\boldphi$, we can compute the $j$th continuous moments
%$\boldm_j$ of $\boldphi$ for $j=0,1,2,\ldots$ by equation ~\eqref{mandM22}.

The eigenvalues of $M_0$ are 1 and $1/2$, and the corresponding eigenvector to the eigenvalue 1 is $[1,0]^T$. Normalization factor is $\sqrt{2}/2$. Hence,
\begin{equation*}
\boldm_0 = \frac{\sqrt{2}}{2} \begin{bmatrix}
                   1  \\
                   0  \\
                \end{bmatrix}
\end{equation*}

We find that
%\begin{equation*}
%d^{j} I_{2} - \left[ M^{+}_{0} + (-1)^j M^{-}_{0} \right]
%= \begin{cases}
%        2^{j}-\frac{1-\sqrt{7}}{4}, & \text{for odd $j$},\\ \noalign{\smallskip}
%        2^{j}-1, \qquad & \text{for even $j$}.
%    \end{cases}
%\end{equation*}
$2^{j} I_{2} - \left[ M^{+}_{0} + (-1)^j M^{-}_{0} \right]$
is invertible for $j=1,2,3,\ldots$.
Hence, we can find the $j$th continuous moments $\boldm_j$ of $\boldphi$ recursively,
for $j=1,2,3,\ldots$.

The $j$th continuous moments $\boldm_j$ of $\boldphi$ for $j=0,1,2$ are
\begin{equation*}
\boldm_0 = \frac{\sqrt{2}}{2} \begin{bmatrix}
                   1  \\
                   0  \\
                \end{bmatrix},\qquad
\boldm_1 = \frac{7\sqrt{2}-\sqrt{14}}{12} \begin{bmatrix}
                   1  \\
                   0  \\
                \end{bmatrix}, \qquad
\boldm_2 = \frac{4\sqrt{2}-\sqrt{14}}{252}\begin{bmatrix}
                   49 \\
                   3\sqrt{3} \\
                \end{bmatrix}.
%\boldm_2 = \begin{bmatrix}
%                   \frac{28\sqrt{2}-7\sqrt{14}}{36}  \\
%                   \frac{4\sqrt{6}-\sqrt{42}}{84}  \\
%                \end{bmatrix}.
\end{equation*}

The discrete moments $N_j$ of $\boldpsi$ for $j=0,1,2$ are
\begin{gather*}
N_0 = \begin{bmatrix}
              0 & 0 \\
              0 & -\frac{\sqrt{3}}{2} \\
            \end{bmatrix}, \qquad
N_1 = \frac{1}{8} \begin{bmatrix}
                       0 & -2\sqrt{7} \\
                       -\sqrt{7} & -7\sqrt{3} \\
            \end{bmatrix}, \qquad
N_2 = \frac{1}{4} \begin{bmatrix}
                       0 & -2\sqrt{7} \\
                       -1+2\sqrt{7} & -7\sqrt{3} \\
                      \end{bmatrix}.
\end{gather*}

%Without computing the positive and negative continuous moments, $\boldn^{+}_j$ and
%$\boldn^{-}_j$ respectively, of $\boldpsi$, we can compute the $j$th continuous
%moments $\boldn_j$ of $\boldpsi$ for $j=0,1,2,\ldots$ by equation ~\eqref{nandN2}.

The $j$th continuous moments $\boldn_j$ of $\boldpsi$ for $j=0,1,2$ are
\begin{gather*}
\boldn_0 = \begin{bmatrix}
              0  \\
              0  \\
            \end{bmatrix}, \qquad
\boldn_1 =  \begin{bmatrix}
                       0 \\
                       0 \\
            \end{bmatrix}, \qquad
\boldn_2 = \frac{1}{168} \begin{bmatrix}
                       0  \\
                       4\sqrt{2}-\sqrt{14} \\
                      \end{bmatrix}.
\end{gather*}

Hence, $\boldpsi$ has 2 vanishing moments and $\boldphi$ provides approximation order 2.
}
$\hfill\Box$
\end{example} % end of example 2

%(I want to plot the two-direction multiscaling functions in examples
%in order to double check the results for $\boldm_0$.)

%\textbf {Acknowledgments.} The author would like to thank Professor Fritz Keinert,
%as well as the Department of Mathematics, at the Iowa State University for their
%hospitality during his visit at the Iowa State University.

%The authors are grateful to the anonymous referees for their comments which greatly improved
%the presentation of this paper.

%\bigskip
%\textbf {Acknowledgment.} The author is grateful to the anonymous
%referees for their comments which improved the presentation
%of this paper.

\bibliography{twodirec}           % calls twodirec.bib
\bibliographystyle{plain}

%\bigskip

%\newpage
%{\bf Experts in the field that are unrelated to the author:}
%\vspace{1cm}
%\address{1. Fritz Keinert, }
%\email{keinert@iastate.edu}

%\vspace{1cm}
%\address{2. Shouzhi Yang, }
%\email{szyang@stu.edu.cn}

%\vspace{1cm}
%\address{3. Youfa Lee, }
%\email{youfalee@hotmail.com}

%\vspace{1cm}
%\address{4. Hong-Tae Shim, }
%\email{hongtae@sunmoon.ac.kr}

%\vspace{1cm}
%\address{5. Dai-Gyoung Kim, }
%\email{dgkim@hanyang.ac.kr}

\end{document}